\documentclass[preprint,12pt]{elsarticle}
\graphicspath{ {./figures/} }
\usepackage{subfig}
\usepackage{xcolor}
\usepackage{amssymb}
\usepackage{amsmath}
\usepackage{hyperref}
\newtheorem{The}{Theorem}[section]
\newtheorem{Lem}{Lemma}[section]
\newtheorem{Ex}{Example}[section]
\newtheorem{Def}{Definition}[section]
\journal{}

\begin{document}
	\begin{frontmatter}
	\title{Stability and Dynamics of Complex Order Fractional Difference Equations}
	\author[1]{Sachin Bhalekar }
	\ead{sachin.math@yahoo.co.in, sachinbhalekar@uohyd.ac.in (Corresponding Author)}
	
	\author[2]{Prashant M. Gade}
	\ead{prashant.m.gade@gmail.com}

	\author[2]{Divya Joshi}
	\ead{divyajoshidj27@gmail.com}

	\address [1]{School of Mathematics and Statistics, University of Hyderabad, Hyderabad, 500046 India}
	\address [2]{Department of Physics, RTM Nagpur University, Nagpur}
\begin{abstract}
We extend the definition of $n$-dimensional difference 
	equations to complex order $\alpha\in \mathbb{C} $.
We investigate the stability of linear systems defined by an $n$-dimensional matrix $A$ and derive conditions for the stability of equilibrium points for
linear systems. For the one-dimensional case where $A =\lambda \in \mathbb {C}$, we find that the stability region,
if any is enclosed by a boundary curve and we obtain a parametric equation for the same. 
Furthermore, we find that 
	there is no stable region 
	if this parametric curve is self-intersecting. 
Even for $ \lambda \in \mathbb{R} $, the solutions can be complex and dynamics in one-dimension is richer than the case for $ \alpha\in \mathbb{R} $.
These results can be extended to $n$-dimensions. For nonlinear systems, we observe that the stability of the linearized system determines the stability of the equilibrium point. 
\end{abstract}
\end{frontmatter}
	\section{Introduction}
The generalization of various mathematical notions such as functions or even operators has  importance that goes beyond mathematical curiosity. The Gamma function is a generalization of factorial for a real argument that is unique under certain constraints \cite{davis1959historical}. It can even be generalized to complex arguments \cite{abramovitz1964handbook}.  It has been used in complex analysis, statistics, number theory, and even string theory in physics. Similarly, the Riemann zeta function was introduced by Euler for real argument and was later extended to complex argument \cite{borwein2008riemann}. It has turned out to be an extremely important function in physics and mathematics \cite{sakhr2003zeta}. The notion of derivatives has even been extended to functional derivatives \cite{bartolotti1982functional}.
 Another example is the q-deformation of numbers and functions \cite{sahoo1993q}. 
 It has found applications in quantum groups and statistical physics \cite{plastino2004liouville}.   
Even q-derivatives and q-integrals have been defined. Of course, these generalizations need not be unique and different generalizations can be used in a different contexts. 

One of the important generalizations of the concept of derivatives has been fractional calculus \cite{kleinz2000child}. 
Several definitions have been proposed 
for extending the concept of fractional derivatives for real numbers. 
Some of them are extended to complex numbers and even
derivatives of fractional complex order have been introduced. Differential equations of complex fractional order (which should be distinguished from complex differential equations) have been studied in the context of viscoelasticity, control systems, etc. Time domain, frequency domain and stability analysis of linear systems represented by differential equations with complex order derivative has been carried out \cite{jacob2016review}. It has found applications such as the design of controller for fractional-order DC motor system \cite{shah2021complex}, in PID controller and low-pass filter \cite{bingi2021design}. Such systems have been found to have large stability regions for certain parameters. The dynamic response of elastic foundations was modeled using complex order differential equations and was useful in predicting the response for various vibration modes over the entire frequency range of interest \cite{makris1994complex}. Particle swarm optimization is a well-studied optimization technique. It has been used in several constraints. Complex order derivatives have found application in this context
as well \cite{pahnehkolaei2021particle}. It has been studied in discrete-time control of linear and nonlinear systems \cite{machado2013optimal}. In biophysics, atrial fibrillation is an important research topic and a mathematical model based on fractional-order complex derivatives has been recently proposed \cite{ugarte2018atrial}. Fractional-order circuit theory has been popular in recent times. It has been extended to fractional-order derivatives in circuit elements \cite{si2017attempt}. Thus complex order differential equations have found applications in several situations as mentioned before and the difference equations of
complex order can be useful in those contexts. 
As Oono and Puri pointed out "Nature gives physicists phenomena, not equations." \cite{oono1988study}, new mathematical tools have always found applications in a variety of fields and these difference equations 
can be useful in numerous phenomena in nature. 
 
Difference equations can be viewed as an attempt to solve the differential equation by finite difference method and the notion of fractional order difference equation has been introduced in this context \cite{atici2010modeling}. 
We note that in fields such as economics and biology, difference equations appear naturally in modeling. Several dynamical phenomena obtained in differential equations are seen in difference equations as well \cite{strogatz}. Many schemes for control of chaos have been applicable to both differential equations as well as maps. The notion of fractional order differential equation has been extended to fractional order difference equation and few definitions have been proposed \cite{deshpande2016chaos}. Dynamics of linear and nonlinear systems have been investigated for fractional-order difference equations \cite{gade2021fractional} and even spatially extended dynamical systems have been defined as well as investigated \cite{pakhare2020emergence}. 

 We are not aware of any attempt to define and study the difference equation of complex fractional order. In this work, we define the difference equation of complex fractional order in Caputo-like definition. We study the stability of linear systems which is an important and useful starting point for understanding the dynamics. We give stability conditions for linear systems and results can be extended to higher dimensions without loss of generality. Finally, we study nonlinear difference equations of complex order and investigate their dynamics.

	\section{Preliminaries}
\begin{Def}
	(see \cite{mozyrska2015transform}). The Z-transform of a sequence $ \{y(n)\}_{n=0}^\infty $ is a complex function given by
	\begin{equation*}
		Y(z)=Z[y](z)=\sum_{k=0}^{\infty} y(k) z^{-k}
	\end{equation*}
	where $z \in \mathbb{C}$ is a complex number for which the series converges absolutely.
\end{Def}
\begin{Def}(see \cite{ferreira2011fractional, bastos2011discrete}).
	Let $ h > 0 ,\; a \in \mathbb{R}$ and $ (h\mathbb{N})_a = \{ a, a+h, a+2h, \ldots\} $.
	For a function $x : (h\mathbb{N})_a \rightarrow  \mathbb{C}$, the forward h-difference operator if defined as 
$$
		(\Delta_h x)(t)=\frac{x(t+h)- x(t)}{h},$$
	 where t	$ \in (h\mathbb{N})_a $.
\end{Def}
	Throughout this article, we take $a = 0$ and $h = 1$. We write $\Delta$ for $\Delta_1 $.
 	Now, we generalize the fractional order operators defined in \cite{mozyrska2015transform,ferreira2011fractional, bastos2011discrete} to include the complex order $\alpha$.

\begin{Def}
	For a function  $x : (h\mathbb{N})_a \rightarrow  \mathbb{C}$ the fractional h-sum of order $\alpha = u +\iota v \in \mathbb{C}, u>0 $ is given by
	\begin{equation*}
		(_{a}\Delta_h^{-\alpha}x)(t) 
		= \frac{h^\alpha}{\Gamma(\alpha)}\sum_{s=0}^{n}\frac{\Gamma(\alpha+n-s)}{\Gamma(n-s+1)} x(a+sh),\\
	\end{equation*}
	where,	$t=a+(\alpha+n)h, \; n \in \mathbb{N_\circ}$.
\end{Def}
	For $h=1$ and $a=0$	, we have
	\begin{eqnarray*}
		(\Delta^{-\alpha}x)(t) 	&=&\frac{1}{\Gamma(\alpha)}\sum_{s=0}^{n}\frac{\Gamma(\alpha+n-s)}{\Gamma(n-s+1)}x(s)\\
			&=&\sum_{s=0}^{n}
		\left(
		\begin{array}{c}
			n-s+\alpha-1\\
			n-s\\
		\end{array}
		\right)
		x(s).
	\end{eqnarray*}                         
	Here, we used the generalized binomial coefficient
	\begin{equation*}
				\left(
			\begin{array}{c}
				\mu \\
				\eta\\
			\end{array}
			\right)
			=\frac{\Gamma(\mu+1)}{\Gamma(\eta+1)\Gamma(\mu-\eta+1)},\\
			\; \mu ,\eta \in \mathbb{C}, \; \text{Re}(\mu)>0,\;\text{and Re}(\eta)>0.
		 	\end{equation*}
			If $n$ $\in$ $\mathbb{N_\circ}$ then
	\begin{eqnarray*}
	   	   \left(
			\begin{array}{c}
			 \mu \\
			  n
			\end{array}
			\right)
			=\frac{(\mu + 1)}{n!\Gamma(\mu-n+1)}
			=\frac{\mu(\mu-1)\ldots(\mu-n-1)}{n!}.
	\end{eqnarray*}
\begin{Def}
	For $n \in \mathbb{N_\circ}$ and $\alpha=u+\iota v \in \mathbb{C}, u>0,$ we define
		\begin{eqnarray*}
		\tilde{\phi}_{\alpha}(n)=
		\left(
		\begin{array}{c}
			n+\alpha-1\\
			n\\
		\end{array}
		\right)
		=(-1)^n
		\left(
		\begin{array}{c}
			-\alpha\\
			n
		\end{array}
		\right).
	\end{eqnarray*}
\end{Def}
\textbf{Note}: The convolution $\tilde{\phi}_{\alpha}*x$ of the sequences $\tilde{\phi}_{\alpha}$ and $x$ is defined as
	\begin{equation*}
		\left(\tilde{\phi}_{\alpha}*x\right)(n)=\sum_{s=0}^{n}\tilde{\phi}_{\alpha}(n-s)x(s)
	\end{equation*}
	\begin{equation*}
		\therefore (\Delta^{-\alpha}x)(n)=(\tilde{\phi}_{\alpha}*x)(n)\\.
	\end{equation*}
\begin{equation*}
		\therefore Z(\Delta^{-\alpha}x)(n)=Z\left(\tilde{\phi}(n)\right)Z(x(n))\\
	=(1-z^{-1})^{-\alpha}X(z),	
\end{equation*}
where $X$ is $Z$ transform of $x$.
\begin{Lem}
	For $\alpha \in \mathbb{C},\; \text{Re}(\alpha)>0$,
	\begin{equation*}
		Z(\tilde{\phi}_{\alpha}(t))=\frac{1}{(1-z^{-1})^{\alpha}}.
	\end{equation*}
\end{Lem}
Proof: We have, 
\begin{eqnarray*}
	Z(\tilde{\phi}_{\alpha}(t))&=&\sum_{j=0}^{\infty}\tilde{\phi}_{\alpha}(j)z^{-j}\\
	&=&\sum_{j=0}^{\infty}\left(
	\begin{array}{c}
		j+\alpha-1\\
		j
	\end{array}
	\right)z^{-j}\\
	&=&\sum_{j=0}^{\infty}(-1)^{j}\left(
	\begin{array}{c}
		-\alpha\\
		j
	\end{array}
	\right)z^{-j}\\
	&=&(1-z^{-1})^{-\alpha}.
\end{eqnarray*}
by using Newton's generalization of Binomial Theorem \cite{niven1969formal,link1}.
\section{Stability Analysis}
We consider the linear fractional order difference equation
\begin{equation}
	(\Delta^{\alpha}x)(t)=(A-I)x(t+\alpha-1), \label{aaa}\\
\end{equation}
where
$	t \in \mathbb{N}_{1-\alpha}=\{1-\alpha,2-\alpha,3-\alpha,\ldots\}, \;
	\alpha \in \mathbb{C},\; \text{Re}(\alpha) \in (0,1), $ 
	$x(t) \in \mathbb{C}^n,\; A$ is $n\times n$ complex matrix and $I$ is  $n\times n$ identity matrix and $x(0)=x_0$.\\
	Initial value problem (\ref{aaa}) is equivalent to \cite{fulai2011existence}        
\begin{equation*}
	x(t)=x_0+\frac{1}{\Gamma(\alpha)}\sum_{s=1+\alpha}^{t+\alpha}\frac{\Gamma(t-s)(A-I)x(s+\alpha-1)}{\Gamma(t-s-\alpha+1)}.
\end{equation*}
Putting $s+\alpha-1=j$, we get
\begin{eqnarray*}
	x(t)&=&x_0+\sum_{j=0}^{t-1}\frac{\Gamma(t-j+\alpha-1)}{\Gamma(\alpha)\Gamma(t-j)}(A-I)x(j)\\
	&=&x_0+(A-I)(\tilde{\phi}_{\alpha}*x)(t-1).
\end{eqnarray*}
\begin{equation}
\therefore \;	x(t+1)=x_0+(A-I)(\tilde{\phi}_{\alpha}*x)(t), \;
	t=0,1,2\ldots.         \label{bbb}
\end{equation}
If $Z(x(t))=X(z),$ then $Z(x(t+1))=zX(z)-zx_0$.\\
Taking Z-transform of (\ref{bbb}), we get
\begin{equation*}
	zX(z)-zx_0
	=\frac{x_0}{1-z^{-1}}+(A-I)\frac{1}{(1-z^{-1})^{\alpha}}X(z)\\.
\end{equation*}
provided, $|z|>1$ \cite{mozyrska2015transform}.
\begin{equation*}
		\therefore[z(1-z^{-1})^{\alpha}I-(A-I)]X(z)\\
	=z(1-z^{-1})^{\alpha-1}x_0,
\end{equation*}
where $|z|>1$. \\
We can solve this equation for $X(z)$ if the matrix 
$(z(1-z^{-1})^{\alpha}I-(A-I))$ is invertible matrix, i.e. if 
$det(z(1-z^{-1})^{\alpha}I-(A-I)) \ne 0$ 
$\forall$ $z$ with $|z|>1$ \cite{elaydi1993stability, desoer2009feedback}. 
Therefore, we have following theorem.
\begin{The}
	The zero solution of (\ref{aaa}) or (\ref{bbb}) is asymptotically stable if and only if all the roots of
$det(z(1-z^{-1})^{\alpha}I-(A-I))=0$
satisfy $|z|<1$.
\end{The}
\subsection{Sketching the boundary of stable region}
Without loss, we can assume that the matrix $(A-I)$ is diagonal.
Suppose that ($\lambda -1$) is an arbitrary entry on the diagonal.
For stability, all the roots of characteristic equation
\begin{equation}
	z(1-z^{-1})^{\alpha}-(\lambda - 1)=0 \label{iii}
\end{equation}
should satisfy $|z|<1$.
On the boundary of stable region, we must have 
$z=e^{\iota t}$, $0\leq t \leq 2\pi$.
Therefore, the characteristic equation (\ref{iii}) becomes 
$ e^{\iota t}(1-e^{-\iota t})^{\alpha}=(\lambda-1)$.
Therefore, we have
\begin{equation*}
	\lambda=2^{\alpha}\left(\sin\frac{t}{2}\right)^{\alpha}e^{\iota\left[\frac{\alpha\pi}{2}+t\left(1-\frac{\alpha}{2}\right)\right]}+1.
\end{equation*}
Therefore, the parametric representation of boundary curve is
\begin{eqnarray}
\gamma(t)=\left(Re[2^{\alpha}\left(\sin\frac{t}{2}\right)^{\alpha}e^{\iota\left[\frac{\alpha\pi}{2}+t\left(1-\frac{\alpha}{2}\right)\right]}]+1,Im[2^{\alpha}\left(\sin\frac{t}{2}\right)^{\alpha}e^{\iota\left[\frac{\alpha\pi}{2}+t\left(1-\frac{\alpha}{2}\right)\right]}]\right),\;t \in [0,2\pi]. \label{eee} 
\end{eqnarray}
If all the eigenvalues of matrix A lie inside this simple closed curve $\gamma(t)$ then the system will be asymptotically stable.
\subsection{Condition for simple curve}
\begin{The}
	The curve $\gamma(t)$ defined by (\ref{iii}) is simple curve for $\alpha=u+\iota v$, $u \in (0,1)$ if and only if $0<v<\sqrt{2u-u^2}$.
\end{The}
\textbf{Proof:} We have
$$ \beta(t)=e^{\iota t}(1-e^{-\iota t})^{\alpha},$$ where $\alpha=u+\iota v$.
There is self intersection or cusp in this parametric curve if and only if $\exists$ $t_1\ne t_2$ such that
\begin{eqnarray*}
	\beta(t_1)&=&\beta(t_2),\; t_1,t_2 \in [0,2\pi].\\
	\iff e^{\iota t_1}(1-e^{\iota t_1})^{u+\iota v}&=&e^{it_2}(1-e^{\iota t_2})^{u+\iota v}\\
	\iff \left(\sin\frac{t_1}{2}\right)^ue^{\frac{vt_1}{2}}e^{\iota [v\log\sin\frac{t_1}{2}+t_1(1-\frac{u}{2})]}&=&\left(\sin\frac{t_2}{2}\right)^ue^{\frac{vt_2}{2}}e^{\iota [v\log\sin\frac{t_2}{2}+t_2(1-\frac{u}{2})]}\\
	\text{and}\\
	v\log\left(\left(\sin\frac{t_1}{2}\right)+t_1\left(1-\frac{u}{2}\right)\right)&=&v\log\left(\left(\sin\frac{t_2}{2}\right)+t_2\left(1-\frac{u}{2}\right)\right)	+2k\pi,\;k \in \mathbb{Z}	\\
	\iff \frac{v^2}{2}&=&\left(1-\frac{u}{2}\right)u-2k\pi u\\
	\therefore v^2&=&(2-u)u-4k\pi u\\
	\iff v&=&\sqrt{(2+4k\pi)u+u^2}
\end{eqnarray*}
This is non real if $k=-1,-2,-3,\ldots$.
Further $v$ is minimum if $k=0$ and we 
have $v=\sqrt{2u-u^2}$ or $(u-1)^2+v^2=1$.
Therefore, if $0<v<\sqrt{2u-u^2}\;$ then $\gamma(t)$ is a simple curve. This completes the proof of theorem.\\
\textbf{Observation}:
\begin{figure}[h]
	\centering
	\includegraphics[scale=1]{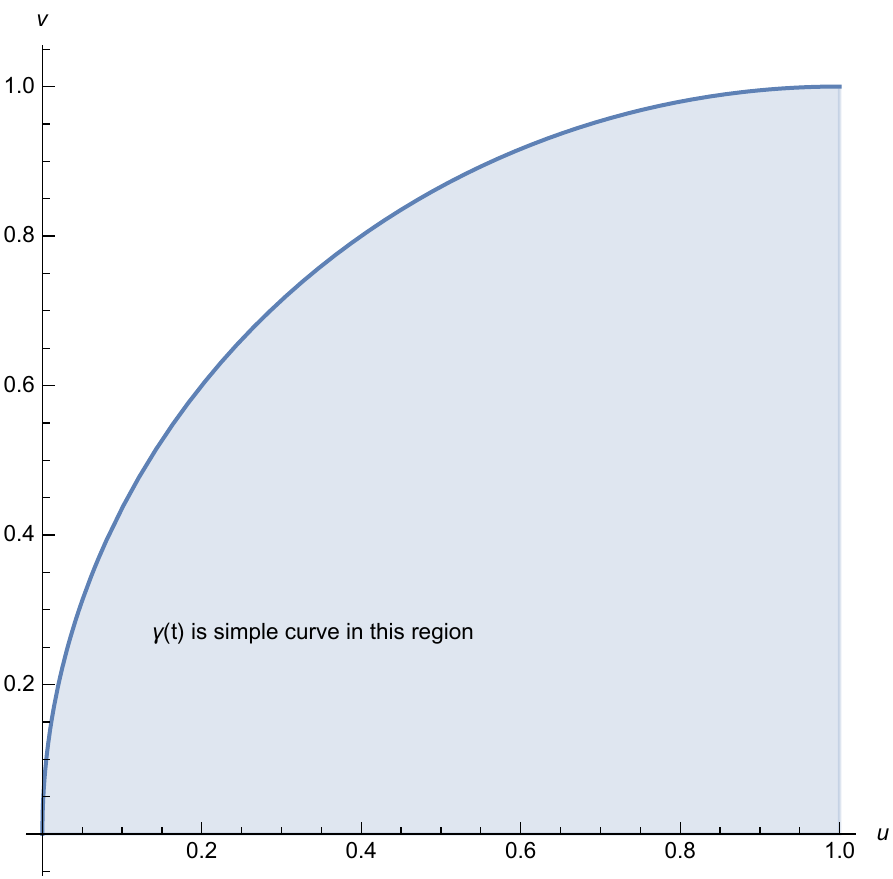}
	\caption{Region for simple curve} \label{fig12}
\end{figure}
\textbf{Observation}:
If there exists multiple points i.e. if $v>\sqrt{2u-u^2}$ then the system is unstable for all eigenvalues.
\section{Illustrative Examples}
In this section we verify the stability results described in the previous section.
\begin{Ex}
	\begin{figure}
		\centering
		\includegraphics[scale=0.6]{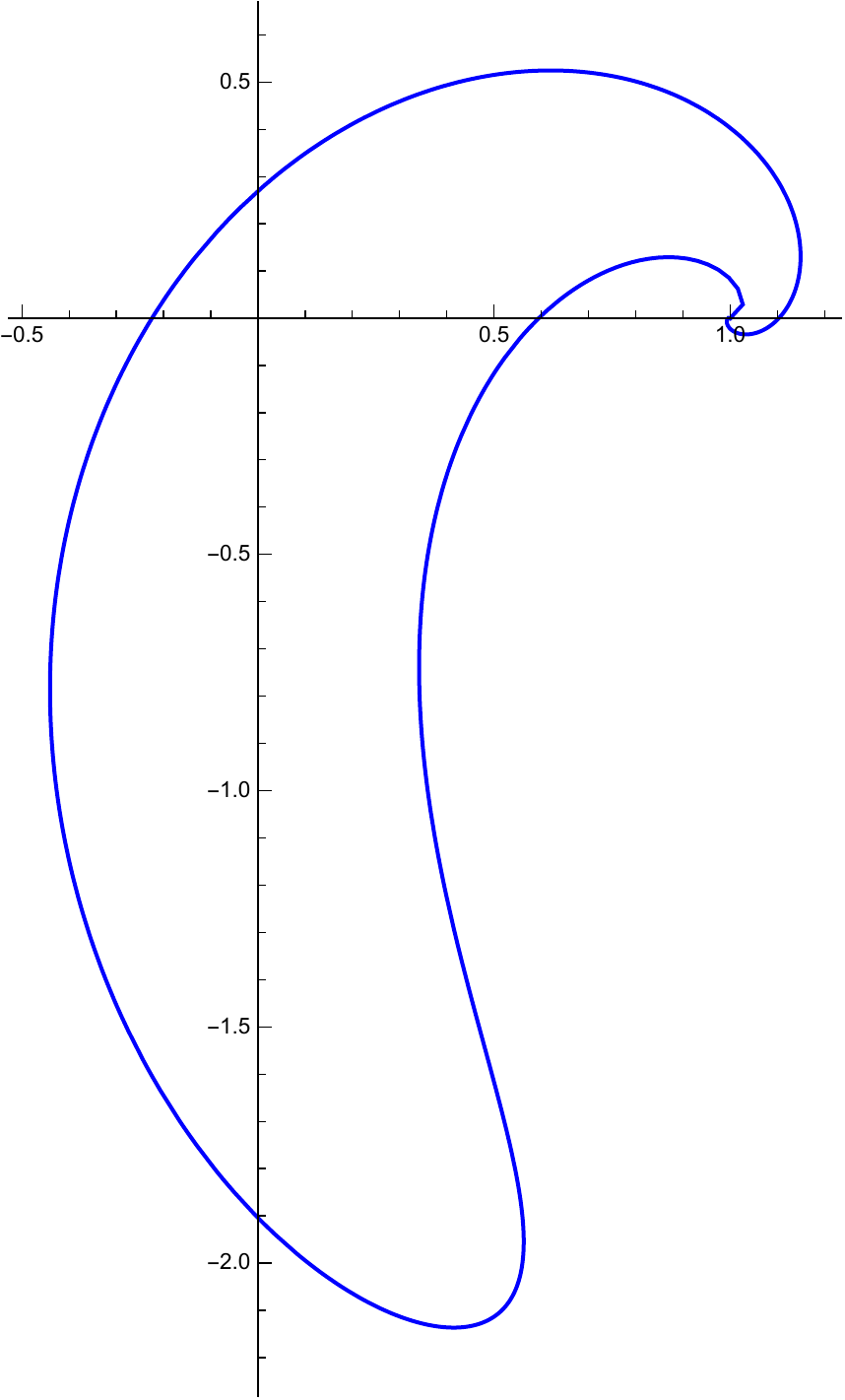}
		\caption{Stable region for $\alpha=e^{\frac{\iota\pi}{4}}$} \label{fig17}
	\end{figure}
	\begin{figure}
			\subfloat{\includegraphics[scale=1]{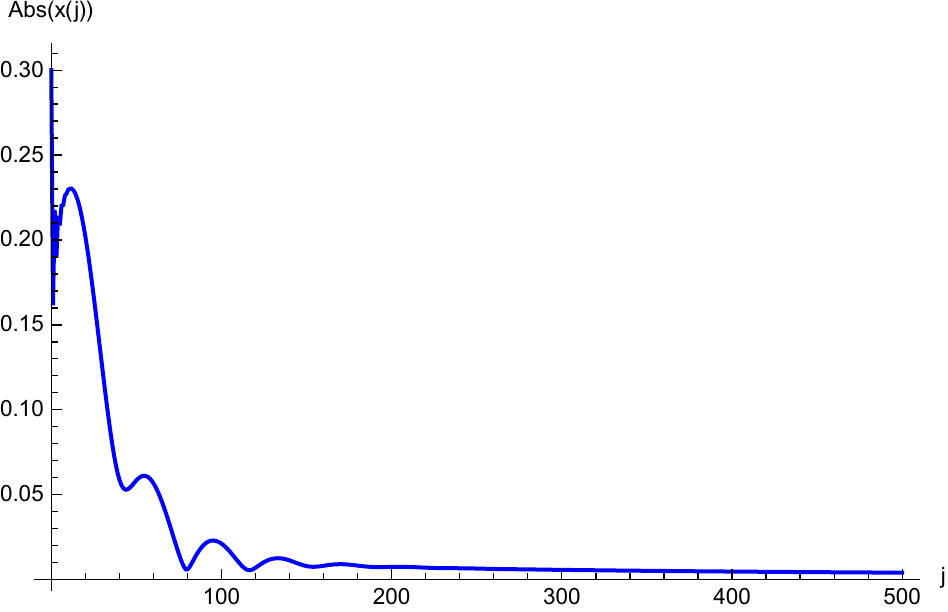}}
			\caption{Example 1 stable solution for $f(x)=(0.2+0.5\iota) x$} 	\label{fig9}
			\subfloat{\includegraphics[scale=1]{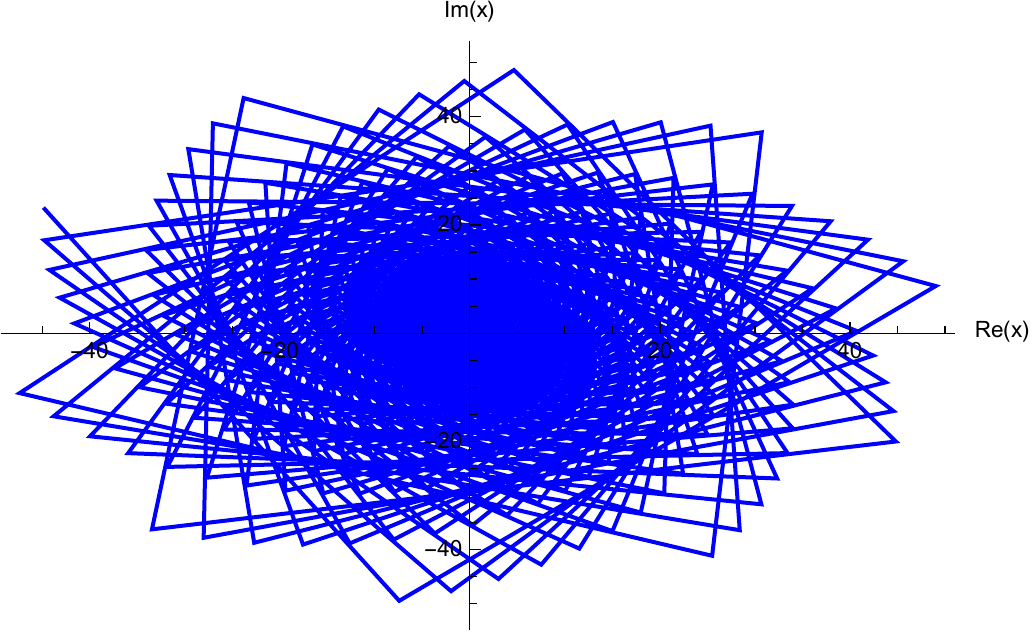}}
			\caption{Example 1 unstable solution for $f(x)=(0.1-2\iota) x$}  \label{fig10}
			\subfloat{\includegraphics[scale=1]{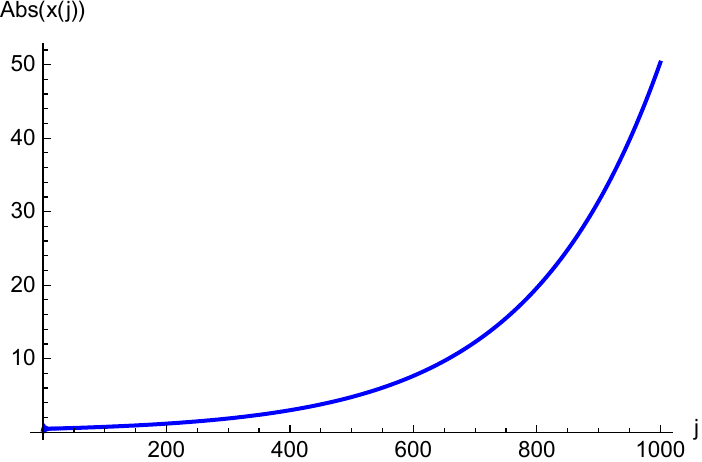}}
			\caption{Example 1 unstable solution for $f(x)=(0.1-2\iota) x$}  \label{fig11}
	\end{figure}
	We take $\alpha=e^{\frac{\iota\pi}{4}}=0.7071(1+\iota)$.
	The stable region for this $\alpha$ is sketched in Figure \ref{fig17}.
Consider the $1$-dimensional system (\ref{bbb}) with $f(x)=(0.2+0.5\iota)x$.
	In this case, $\lambda=0.2+0.5\iota$ 
	is inside the stable region of this system as shown 
	in Figure \ref{fig9}.
	On the other hand, the eigenvalue of $\lambda=0.1-2\iota)x$ lies outside the stable region. The unstable trajectory of this system is traced in Figures \ref{fig10} and \ref{fig11}.
\end{Ex}
\begin{Ex}
	\begin{figure}[h]
		\centering
		\includegraphics[scale=0.9]{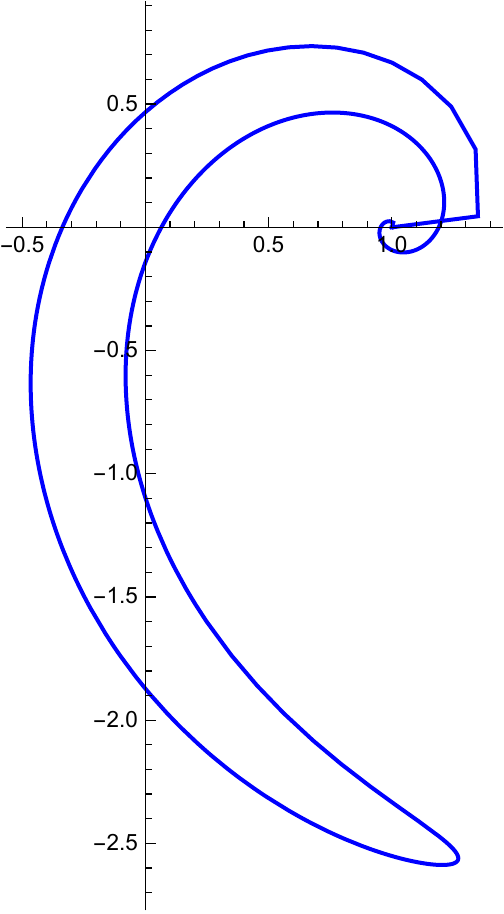}
		\caption{Multiple curve for $\alpha=0.4+0.9\iota$} \label{fig13}
	\end{figure}
	\begin{figure}
		\subfloat{\includegraphics[scale=1]{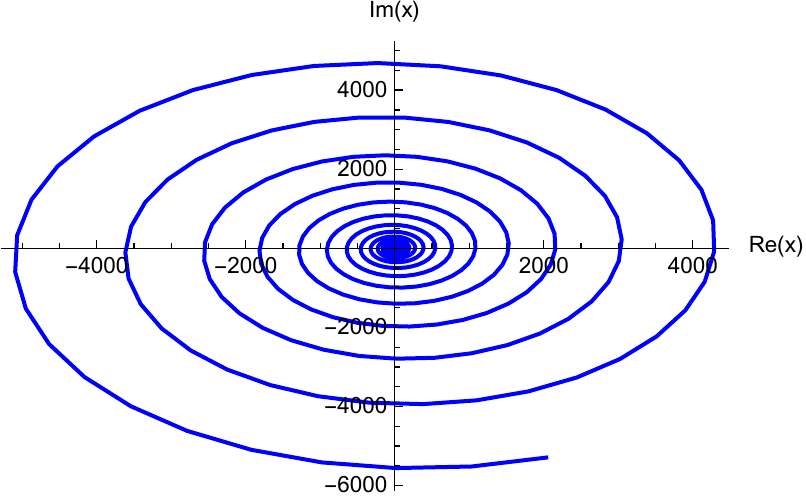}}
		\caption{Example 2 unstable trajectory for $\lambda=1.1-0.1\iota$} \label{fig14}
		\subfloat{\includegraphics[scale=0.9]{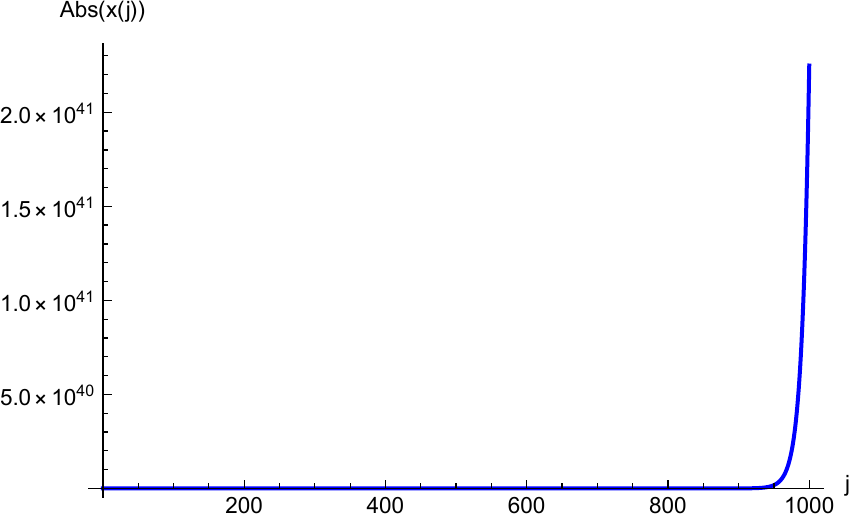}}
		\caption{Example 2 unstable solution for $\lambda=1.1-0.1\iota$}\label{fig15}
		\subfloat
		{\includegraphics[scale=1]{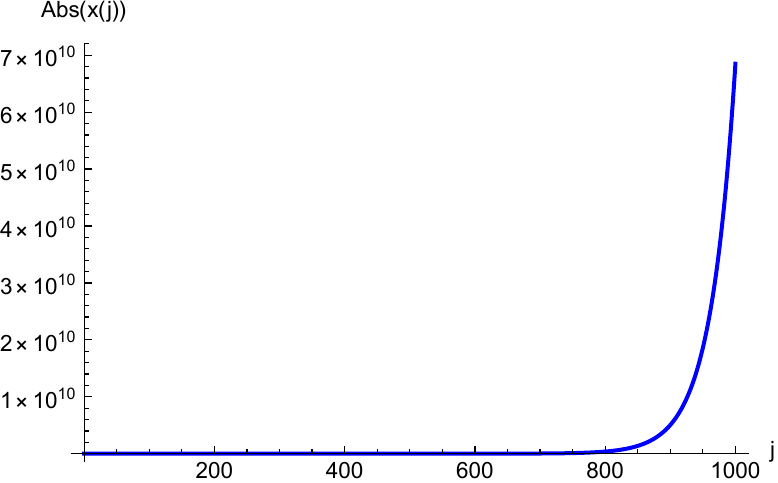}}
		\caption{Example 2 unstable solution for $\lambda=0.5+0.1\iota$}\label{fig16}
	\end{figure}
	Now we take $\alpha=0.4+0.9\iota$.
Here $v>\sqrt{2u-u^2}$.
Therefore, the boundary 
curve $\gamma(t)$ will have multiple points as shown in Figure \ref{fig13}.
	In this case, we observe the unstable solutions for all eigenvalues. 
We take $f(x)=\lambda x$ with $\lambda=1.1-0.1\iota,-0.1-0.1\iota$ 
	and $0.5+0.1\iota$ inside and outside the curve $\gamma(t)$. The unstable solutions for these values are sketched in Figures \ref{fig14}, \ref{fig15} and \ref{fig16} respectively.
\end{Ex}
\section{Nonlinear System}
We consider, 
\begin{equation}
	x(t)=x_0+\frac{1}{\Gamma(\alpha)}\sum_{j=0}^{t-1}\frac{\Gamma(t-j+\alpha-1)}{\Gamma(t-j)}[f(x(j))-x(j)],\; t=1,2,\ldots.   \label{fff}
\end{equation}
where $x(t) \in \mathbb{C}^n$ and $f:\mathbb{C}^n \rightarrow \mathbb{C}^n$ is continuously differentiable.
	A steady state solution $x_*$ of (\ref{fff}) is a complex number 
	satisfying $f(x_*)=x_*$.

Following definitions are generalizations of those given in \cite{elaydi2006introduction, hirsch2012differential}.
\begin{Def}
We say that $x_*$ is stable if for each $\epsilon>0,\; \exists\; \delta>0$ such that $\Vert x_0-x_*\Vert <\delta$  $\implies\; \Vert x(t)-x_*\Vert<\epsilon$, $t=1,2,\ldots$
\end{Def}
\begin{Def}
	Equilibrium point $x_*$ is asymptotically stable if it is stable and $\exists\; \delta>0$ such that $\Vert x_0-x_*\Vert<\delta \; \implies \lim_{t\rightarrow \infty}x(t)=x_*$.
\end{Def}
\textbf{Note}: If $x_*$ is not stable then it is unstable. \\
The linearization of nonlinear system (\ref{fff}) in the neighborhood of $x_*$ is given by
\begin{equation}
	x(t)=x_0+\sum_{j=0}^{t-1}\frac{\Gamma(t-j+\alpha-1)}{\Gamma(\alpha)\Gamma(t-j)}(A-I)x(j)  \label{ggg}
\end{equation} 
where $A=f'(x_*)$ is the Jacobian matrix.

The local stability properties of equilibrium point $x_*$ of (\ref{fff}) are same as those of linearization (\ref{ggg}) i.e. the equilibrium point $x_*$ is asymptotically stable if all the eigenvalues of Jacobian $A$ lie inside the stable region.
\subsection{Complex order logistic map}
We consider 
$$f(x)=\lambda x(1-x),$$ where $\lambda>0$.\\
The logistic map of complex order $\alpha$, with Re$(\alpha)>0$ is given by
\begin{equation}
	x(t)=x_0+\sum_{j=0}^{t-1}\frac{\Gamma(t-j-\alpha-1)}{\Gamma(\alpha)\Gamma(t-j)}[f(x(j))-x(j)],\; t=1,2,\ldots. \label{ccc}
\end{equation}   
Equilibrium points of (\ref{ccc}) are $x_{1*}=0$ and $x_{2*}=(\frac{\lambda-1}{\lambda})$.
Linearization of (\ref{ccc}) in the neighborhood of equilibrium $x_*$ is given by
\begin{equation}
	x(t)=x_0+\sum_{j=0}^{t-1}\frac{\Gamma(t-j-\alpha-1)}{\Gamma(\alpha)\Gamma(t-j)}[(f'(x_*)-1)x(j)]. \label{ddd}
\end{equation}
\textbf{Stability of $x_{1*}=0$:}\\
Here $f'(x_{1*})=f'(0)=\lambda$.
Therefore, the equilibrium point $x_{1*}$ is 
asymptotically stable if $\lambda$ 
is inside the stable region bounded by (\ref{eee}).\\
\textbf{Stability of $x_{2*}=0$:}\\
In this case,
\begin{eqnarray*}
	f'(x_{2*})&=&\lambda-2\lambda\left(\frac{\lambda-1}{\lambda}\right)\\
	&=&2-\lambda.
\end{eqnarray*}
Therefore, $x_{2*}$ is asymptotically stable if ($2-\lambda$) lies inside the stable region bounded by (\ref{eee}).\\
\begin{figure}[h]
	\centering
	\includegraphics[scale=0.8]{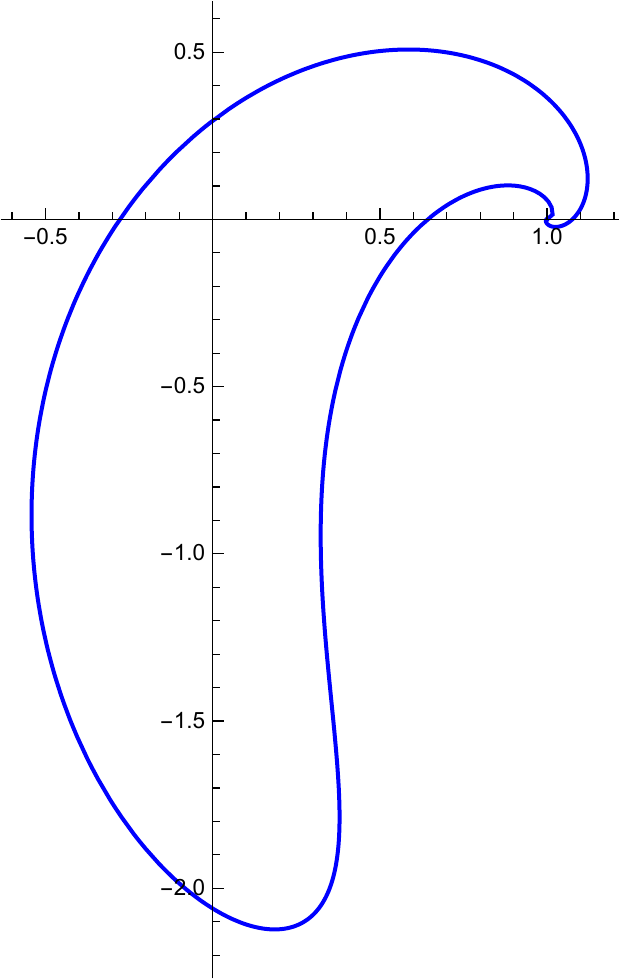}
	\caption{Stability region for logistic map} \label{fig1}
\end{figure}
\begin{figure}[h]
	\subfloat
	{\includegraphics[scale=1]{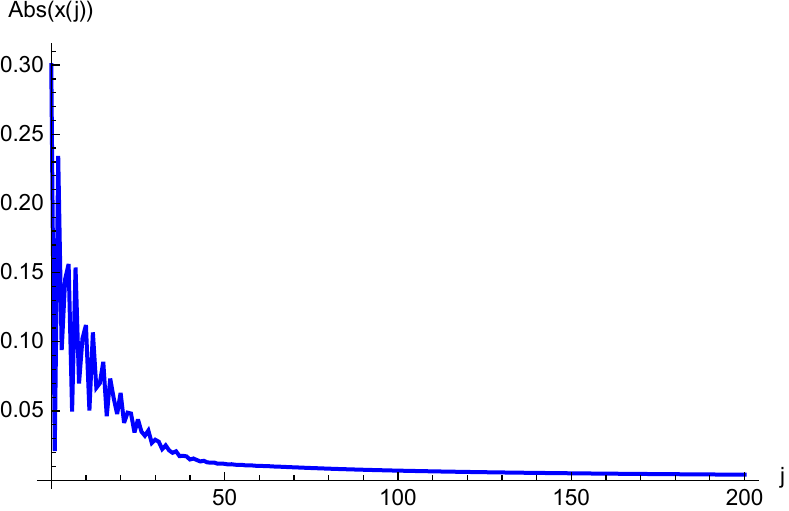}}
	\caption{$x_{1*}$ is asymptotically stable for $\lambda=-0.1$}\label{fig2}
	\subfloat{\includegraphics[scale=1]{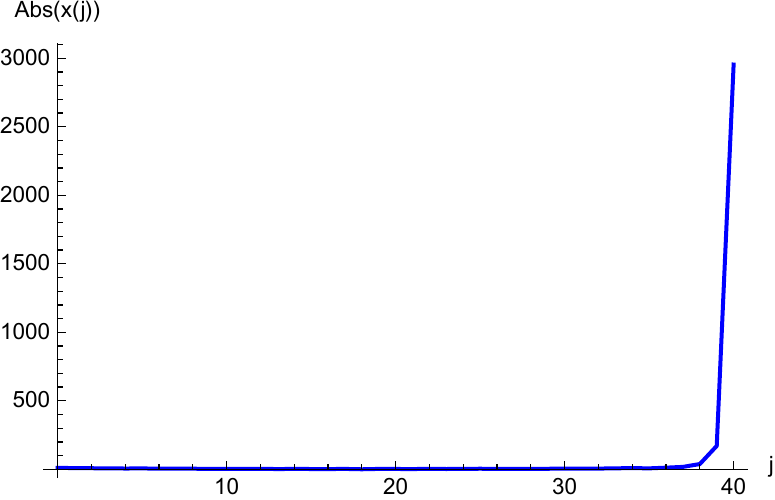}}
	\caption{$x_{2*}$ is unstable for $\lambda=-0.1$}\label{fig3}
	\subfloat{\includegraphics[scale=1]{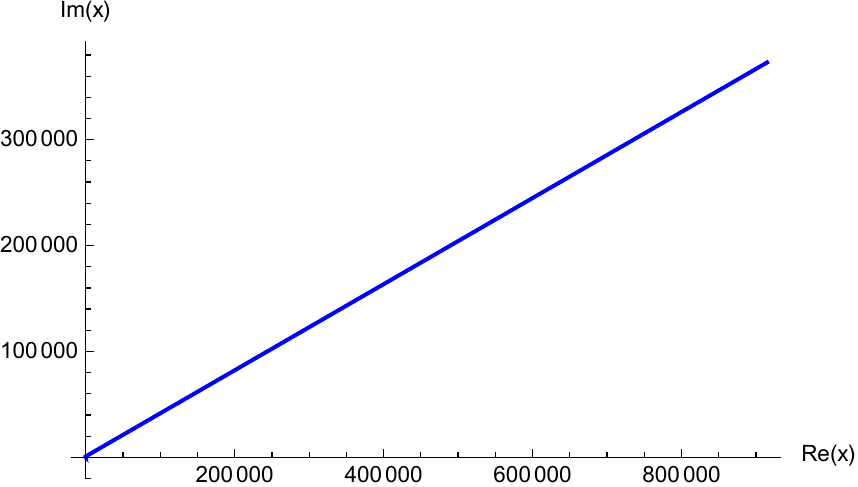}}
	\caption{$x_{1*}$ is unstable for $\lambda=1.5$}\label{fig4}
	\subfloat{\includegraphics[scale=1]{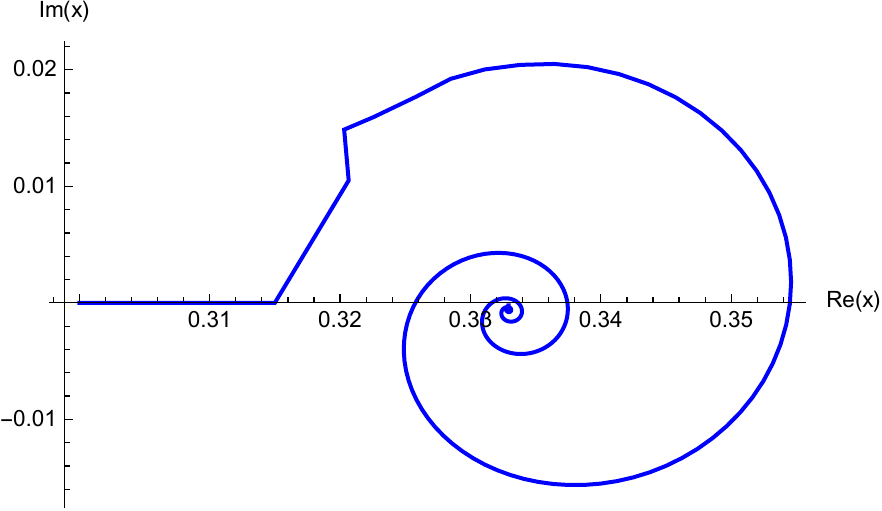}}
	\caption{$x_{2*}$ is asymptotically stable for $\lambda=1.5$}\label{fig5}
\end{figure}
Let us take $\alpha=0.8+0.7\iota $.
The stable region is given in Figure \ref{fig1}.
Therefore, $x_{1*}$ is asymptotically stable if 
$\lambda \in (-0.2774,0.6432)\cup (1,1.0754)$ and 
the $x_{2*}$ is asymptotically stable 
if $\lambda \in (0.9246,1)\cup (1.3568,2.2774)$.
For $\lambda=-0.1$, $x_{1*}$ is asymptotically stable whereas $x_{2*}=11$ is unstable. In Figure \ref{fig2}, we take $x_0=0.3$ where the trajectory converges to $x_{1*}$. We take $x_0=10.2$ in Figure \ref{fig3} and show that the trajectory is diverging. The equilibrium point $x_{1*}$ is unstable for $\lambda=1.5$ whereas $x_{2*}=0.3333$ is asymptotically stable. This result is verified in Figures \ref{fig4} and \ref{fig5} by selecting appropriate initial conditions viz. $x_0=0.3$ and $x_0=-0.1$ respectively.
\subsection{Two-dimensional system}
Now, we consider two-dimensional system
\begin{eqnarray}
		x(t)=x_0+\frac{1}{\Gamma(\alpha)}\sum_{j=0}^{t-1}\frac{\Gamma(t-j+\alpha-1)}{\Gamma(t-j)}[f_{1}(x(j),y(j))-x(j)], \nonumber \\
		y(t)=y_0+\frac{1}{\Gamma(\alpha)}\sum_{j=0}^{t-1}\frac{\Gamma(t-j+\alpha-1)}{\Gamma(t-j)}[f_{2}(x(j),y(j))-y(j)],\label{jjj}\\ 
		 t=1,2,\ldots \nonumber
\end{eqnarray}
where,  
$	f_{1}(x,y)=\lambda x(y+1)+\mu(x^2+1)y$ and $
	f_{2}(x,y)=\lambda y(x+1)-\mu(y+1)^2x.$
Clearly, origin $(0,0)$ is an equilibrium point.
The Jacobian matrix at origin
\begin{eqnarray*}
	J=\left[
	\begin{array}{cc}
		\lambda&\mu\\
		\-\mu&\lambda
	\end{array}\right]
\end{eqnarray*}
\begin{figure}[h]
	\centering
	\includegraphics[scale=0.9]{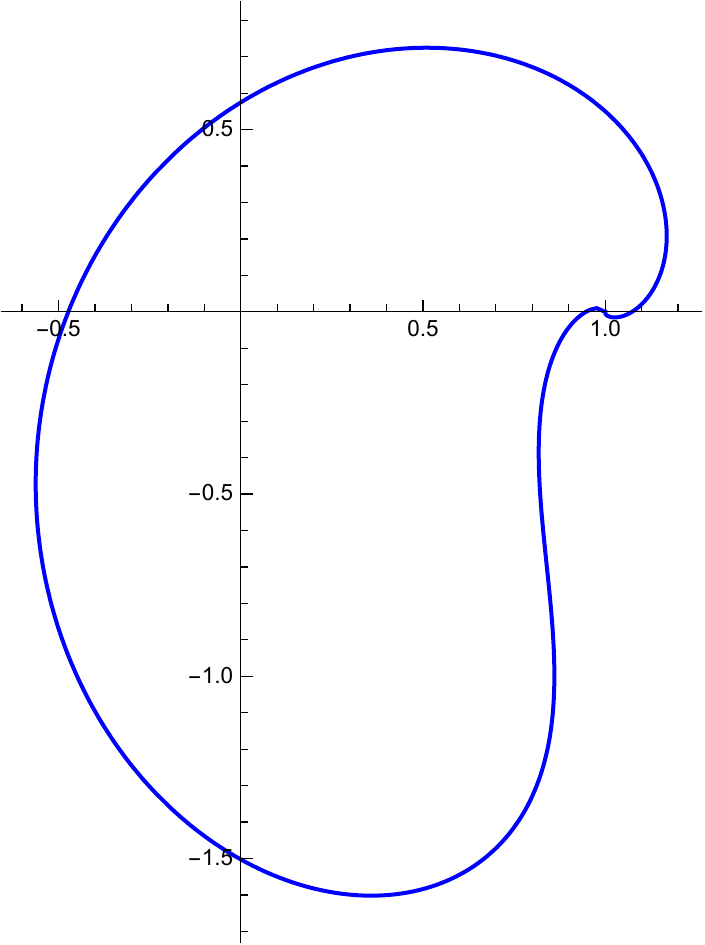}
	\caption{Stable region for $2$-dimensional system (\ref{jjj})}
	\label{fig6}
\end{figure}
\begin{figure}
	\subfloat{\includegraphics[scale=0.9]{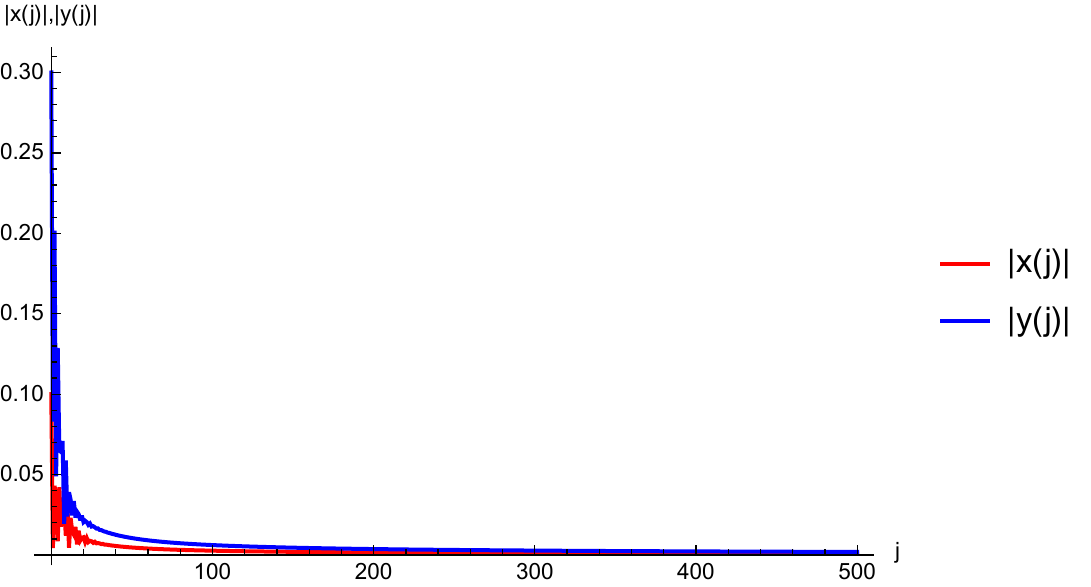}}
	\caption{Stable orbits for $(\lambda,\mu)=(-0.2,0.1)$}
	\label{fig7}
	\subfloat{\includegraphics[scale=0.9]{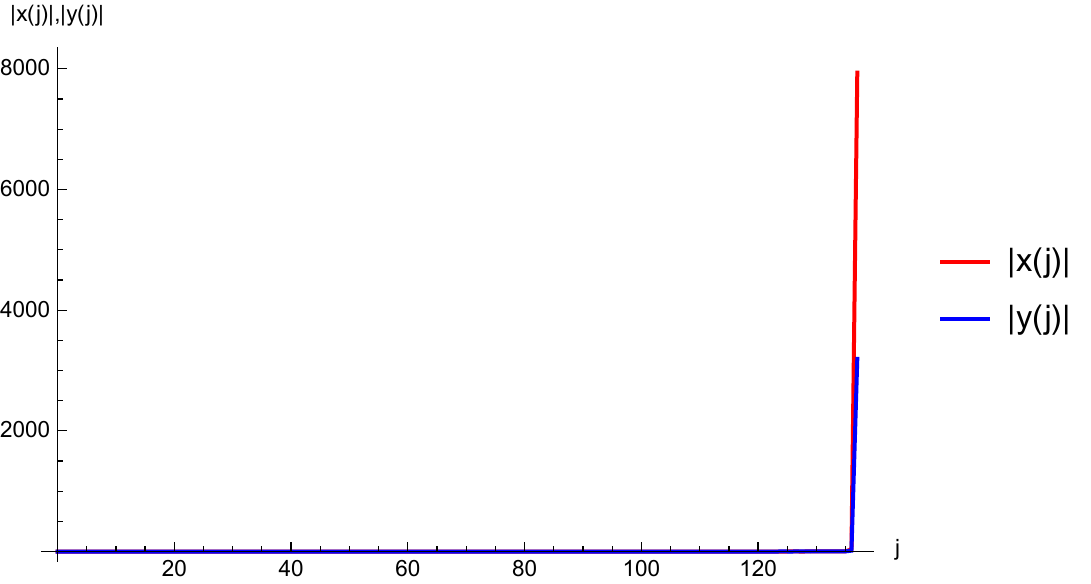}}
	\caption{Unstable orbits for $(\lambda,\mu)=(-0.2,0.5)$}
	\label{fig8}
\end{figure}
has eigenvalues $\lambda \pm \iota \mu$.
For $\alpha=0.7+0.4\iota$, the stable region is shown in Figure \ref{fig6}.
For the value $(\lambda,\mu)=(-0.2,0.1)$, both the eigenvalues 
$-0.2\pm0.1\iota$ of $J$ lie inside the stable region. 
The stable orbits $|x(j)|$ and $|y(j)|$ are shown in Figure \ref{fig7}.
If we set $(\lambda,\mu)=(-0.2,-0.5)$ then the eigenvalue $-0.2-0.5\iota$ lies inside the stable region whereas the eigenvalue $-0.2+0.5\iota
$ lies outside the stable region.
Therefore, the equilibrium is unstable.
The unstable orbits are shown in Figure \ref{fig8}.
\section{Results and Conclusion}
We have used Z-transform for carrying out stability analysis of equilibrium points. This technique is fairly general. Stability for the linear system is achieved if the zero solutions of the system are asymptotically stable i.e., the roots of equation (\ref{iii}) satisfy $|z|<1$. If all the roots of the system lie within the boundary curve defined by (\ref{eee}) then it is considered to be asymptotically stable otherwise, unstable. We have demonstrated that this result can be extended to non-linear systems using the logistic map. For a non-linear system, the first step is to linearize it around its equilibrium points. These points are asymptotically stable if the eigenvalues of the Jacobian matrix lie inside the stable region. For higher dimensions, we considered the example of a two-dimensional system and found similar results for stability analysis. All the eigenvalues of the Jacobian matrix must lie inside the curve for the equilibrium point to be stable as seen in the examples. This criterion is very similar to one obtained for integer-order difference equations though the stability region is very different.

However, we note that the qualitative dynamics is dissimilar in many contexts. For example, for $x_{n+1}=\lambda x_n$, we observe monotonic decrease or increase for $\lambda \in \mathbb{R}$. On the other hand, trajectories can spiral in or out even for $\lambda \in \mathbb{R}$ for complex order difference equations even for real initial conditions. It is possible that these difference equations can lead to stable limit cycles or quasi-periodic cycles in nonlinear systems even in $1$-dimension and it is of interest in the future if it is possible in one or higher dimensions. It has been shown that periodic orbits do not exist for fractional-order differential equations for real order \cite{saleh2012simplification}. Complex order differential has been studied in several contexts as mentioned in the introduction. In viscoelastic systems, some researchers have demanded that the output should be real for real input, and hence the dynamics is governed by the sum of two complex order differential equations of complex conjugate order \cite{atanackovic2016complex}. Our work can be easily extended in this direction and the stable region, in this case, would be the intersection set of the two difference equations.

This work is mainly motivated by mathematical curiosity and we have generalized the notion of difference equations to complex order. The limited studies in nonlinear systems do not show chaos or quasi-periodic cycles. It would be interesting to further examine whether the system shows chaos under certain conditions.

\section{Acknowledgment}
P. M. Gade thanks DST-SERB for financial assistance (Ref. EMR/2016/006686 and CRG/2020/003993).
S. Bhalekar acknowledges the Science and Engineering Research Board (SERB), New Delhi, India for the Research Grant (Ref. MTR/2017/000068) under Mathematical Research Impact Centric Support (MATRICS) Scheme and the University of Hyderabad for Institute of Eminence-Professional Development Fund (IoE-PDF) by MHRD (F11/9/2019-U3(A)).

\bibliography{reference.bib}
\bibliographystyle{elsarticle-num}

\end{document}